\documentclass[11pt]{article}

\usepackage{amsmath,wasysym,cite}
\usepackage{amssymb}
\usepackage{eucal,xcolor}

\begin{document}

\baselineskip 16pt

\title{Finite groups with generalized Ore supplement
conditions for primary subgroups}

\author{Wenbin Guo\thanks{Research of the first author is supported by
a NNSF grant of China (Grant \# 11371335) and Wu Wen-Tsun Key Laboratory of Mathematics, USTC, Chinese Academy of Sciences.}\\
{\small School of Mathematical Sciences, University of Science and
Technology of China,}\\ {\small Hefei 230026, P. R. China}\\
{\small E-mail:
wbguo@ustc.edu.cn}\\ \\
{ Alexander  N. Skiba \thanks{Research of the second author
 supported by Chinese Academy of Sciences Visiting Professorship for Senior
 International Scientists (grant No. 2010T2J12).}}\\
{\small Department of Mathematics,  Francisk Skorina Gomel State University,}\\
{\small Gomel 246019, Belarus}\\
{\small E-mail: alexander.skiba49@gmail.com}}

\date{}

\maketitle

\begin{abstract}  We consider some  applications  of the theory of
generalized Ore supplement conditions in the study of  finite groups.

\end{abstract}

\let\thefootnoteorig\thefootnote
\renewcommand{\thefootnote}{\empty}


\let\thefootnote\thefootnoteorig

\renewcommand{\thefootnote}{\empty}

\section{Introduction}

Throughout this paper, all groups are finite and $G$ always denotes
a finite group. Moreover, $p$ is always supposed to be a prime.
A subgroup  $H$ of $G$ is said to be  \emph{$S$-quasinormal}  in $G$ if $H$ permutes with each Sylow
subgroup $P$ of $G$, that is, $HP=PH$. We use $\frak{U}$  to denote the  class  of all supersoluble
groups.\par

Let $\frak{F}$  be a class of groups.
A chief factor $H/K$ of $G$ is said to be  $\frak{F}$-central in $G$
if  $(H/K)\rtimes (G/C_{G}(H/K))\in \frak{F}$. The product of all
normal subgroups $N$ of $G$ such that every chief factor of $G$
below $N$ is $\frak{F}$-central in $G$ denoted by $Z_\frak{F}(G)$
and called the $\frak{F}$-hypercentre of $G$. By the
Barnes-Kegel  theorem \cite[IV, 1.5]{DH}, for any group
$G\in\frak{F}$ we have $Z_\frak{F}(G)=G$ provided $\frak{F}$ is a formation.

If $G=HB$, then $B$ is said to be a \emph{supplement} of $H$ in $G$.
Since $HG=G$, it makes sense   to consider only the supplements $B$
with some restrictions on $B$. For example, we often deal with  the
situation when for a supplement $B$ of $H$ we have $H\cap B=1$. In
this case, $B$ is said to be a \emph{complement } of $H$ in $G$ and
$H$ is said to be \emph{complemented } in $G$; if, in addition,
 $B$ is also  normal in $G$, then $B$ is said to be a \emph{ normal
complement}  of $H$ in $G$. Note that if  $H$ either is normal in
$G$ or has a normal complement in $G$, then, clearly, $H$ satisfies
the following: there are normal subgroups $T$ and $S$ of $G$ such
that $G=HT$, $S\leq H$ and $H\cap T\leq S$ (see Ore \cite{Ore}); in
this case we say $H$ satisfies the \emph{Ore supplement condition}
in $G$. It is clear that $H$ satisfies the Ore supplement
condition if and only if $H/H_{G}$ has a normal complement in
$G/H_{G}$. In the paper \cite{Wang}, the subgroups which satisfy the
Ore supplement condition were called \emph{$c$-normal}.

It was discovered  that many important for applications classes of
groups (for example, the classes of all soluble, supersoluble,
nilpotent, $p$-nilpotent, metanilpotent, dispersive in the sense
of Ore \cite{Ore} groups (i.e., groups having a Sylow tower)) may be characterized in the terms of the
Ore supplement condition or in the terms of some generalized Ore
supplement conditions. It was the main motivation for introducing,
studying and applying the generalized Ore supplement
conditions of various type.  But, in fact, all recent results in this line of researches are based on the
ideas of the papers \cite{Bal5, GuoManu1, Ski}, in which the authors analyze three  fundamentally
different generalizations of the Ore supplement condition.  A subgroup $H$
of $G$ is called:  \emph{$c$-supplemented} in $G$ \cite{Bal5}
provided $H/H_{G}$ is complemented in $G/H_{G}$;
\emph{${\frak{F}}$-supplemented} \cite{GuoManu1} in $G$ provided there is a
supplement  $T/H_{G}$ of $H/H_{G}$  in $G/H_{G}$ such that
$(H/H_{G})\cap (T/H_{G}) \leq Z_{\frak{F}}(G/H_{G})$;  \emph{weakly $S$-permutable} in
$G$ \cite{Ski} provided there is a subnormal subgroup $T$ of $G$ such that
$HT=G$ and $H\cap T\leq S\leq H$ for some $S$-quasinormal subgroup
$S$ of $G$. Finally, also  we often meet the situation when a
subgroup $H$ has a supplement $T$ in $G$ such that $T\in \frak{F}$.
In spite of the four supplement conditions are  quite
different, there are a lot of similar results, in which we meet one
of these ones.

It is known for example that a soluble group $G$ is supersoluble provided $G$
 has a normal subgroup $E$ with supersoluble quotient $G/E$ such that
 for every
maximal subgroup $H$ of every  Sylow subgroup of $F(E)$ at least one of
the following holds:

(I) $H$ is a $CAP$-subgroup of $G$ \cite{Ezq1}, that is, $H$ either covers or
avoids each chief factor of $G$  (see
\cite[p.~37]{DH});

(II) $H$ is complemented \cite{B-G} or, at least, $c$-complemented in $G$  \cite{Bal5};

(III)  $H$ has a supersoluble supplement in $G$  \cite{GSSII};

(IV)  $H$ is  ${\frak{U}}$-supplemented in $G$ \cite{GuoManu1};

(V) $H$ is  weakly $S$-permutable in $G$ \cite{Ski};

(VI) $H$ is a modular element (in the sense of Kurosh \cite[p.~43]{Schm})  of the
lattice of all subgroups of $G$ \cite{Ukr1}.

The similarity of these results, as well as the similarity  of many other  results of this kind,   makes
 natural to ask:

{\bf Question A.}  {\sl Is there a condition which generalizes
all these conditions on the maximal subgroups of
the Sylow subgroups, and under which $G$ is still supersoluble?}

In fact, the solution of this problem is based on the
concept of the $\frak{F}$-hypercentre and on the idea of the
subgroup functor.

{\bf  Definition 1.1.} Let $K\leq H $ be   subgroups of $G$. Then we
say that that the  pair $(K, H)$ satisfies the
\emph{$\frak{F}$-supplement condition}  in $G$ if $G$ has a subgroup
$T$ such that $HT=G$ and $ H\cap T \subseteq KZ_\frak{F}(T)$.

Recall that a \emph{subgroup functor} is a function $\tau$ which
assigns to each group $G$ a set ${\tau}(G)$  subgroups
 (perhaps  consisting of a single element) of $G$
satisfying $1\in {\tau}(G)$ and  ${\theta}
({\tau}(G))={\tau}({\theta} (G))$ for any isomorphism ${\theta}:G\to
G^{*}$. If $H\in {\tau}(G)$, then we say that $H$ is a
${\tau}$-\emph{subgroup} of the group $G$.

For our goal, we need  the following realization of Definition 1.1.

{\bf  Definition 1.2.}   Let $\bar{G}=G/H_{G}$ and
$\bar{H}=H/H_{G}$, where $H$ is a subgroup of $G$. Let $\tau  $ be a
subgroup functor. Then we say that that $H$ is   \emph{
$\frak{F}_{\tau}$-supplemented} in $G$ if  for some $\tau $-subgroup
$\bar{S}$ of $\bar{G}$ contained in $\bar{H}$ the pair $(\bar{S},
\bar{H})$  satisfies the $\frak{F}$-supplement condition    in
$\bar{G}$.

We show that this  concept gives the positive answer to the first
part of Question A. First note that if $\tau(G)$ is the set of $CAP$-subgroups of $G$ and
$H\in \tau(G)$, then $H/H_{G}$ is a $CAP$-subgroup of $G/H_{G}$ (see
Example 1.4 below), so the triple  $(\bar{H}, \bar{S}, \bar{T})$, where
$\bar{S}=\bar{H}  $ and $ \bar{T}=\bar{G}$,
satisfies Definition 1.2. Clearly,  a subgroup $H$ is $c$-supplemented in
$G$ if and only if it is $\frak{F}_{\tau}$-supplemented in $G$,
where $\frak{F}=(1)$ is the class of all identity groups and
$\tau(G)$ is the set of all normal subgroups of $G$. If $H$ has a
supersoluble supplement $T$ in $G$, then $(H/H_{G})\cap
(H_{G}T/H_{G})=(H\cap T)H_{G}/H_{G}\leq TH_{G}/H_{G}\simeq
T/H_{G}\cap T$, so $H$ is $\frak{U}_{\tau}$-supplemented in $G$ for
any subgroup functor $\tau$. If $H$ is $\frak{U}$-supplemented in
$G$, then similarly we get that $H$ is
$\frak{U}_{\tau}$-supplemented in $G$ for any subgroup functor
$\tau$.  Finally, if $H$ is weakly $S$-permutable in $G$, and $T$ is
a subnormal subgroup of $G$ such that $HT=G$ and $H\cap T\leq S\leq
H$ for some $S$-quasinormal subgroup $S$ of $G$, then $(H/H_{G})\cap
(TH_{G}/H_{G})=(H\cap T)H_{G}/H_{G}\leq  SH_{G}/H_{G}$, where
$SH_{G}/H_{G}$ is a $S$-quasinormal subgroup of $G/H_{G}$ (see
Example 1.6 below). Hence $H$ is $\frak{F}_{\tau}$-supplemented in
$G$ for any formation $\frak{F}$ and for the  subgroup functor
$\tau$, which assigns to each group $G$ the  set ${\tau}(G)$ of all
$S$-quasinormal subgroup of $G$.

Our  next goal is to give the positive answer to the second part of
Question A.

But first,  we define  some   subgroup functors which will be
used in applications of the results.

{\bf   Definition 1.3.} Let ${\tau}$ be a subgroup functor. Then we
say that ${\tau}$ is:

(1) \emph{ Inductive}  provided $HN/N\in {\tau} (G/N)$ whenever
$H\in {\tau} (G)$ and $N\trianglelefteq G$.

(2) \emph{Hereditary} provided $\tau$  is inductive and $H\in \tau
(E)$ whenever $H\leq E\leq G$ and $H\in \tau (G)$.

(3) \emph{$\Phi$-regular} (respectively \emph{$\Phi$-quasiregular})
provided for any primitive group $G$, whenever $H\in {\tau} (G)$ is
a $p$-group and $N$ is a minimal normal (minimal normal abelian, respectively)
subgroup of $G$, then $|G:N_{G}(H\cap N)|$ is a power of $p$.

(4) \emph{Regular} or  a \emph{Li-subgroup functor} \cite{Li}
provided for any group $G$,  whenever $H\in {\tau} (G)$ is a
$p$-group and $N$ is a minimal normal subgroup of $G$, then
$|G:N_{G}(H\cap N)|$ is a power of $p$.

(5)   \emph{Quasiregular } provided for any group $G$, whenever
$H\in {\tau} (G)$ is a $p$-group and $N$ is an abelian minimal
normal subgroup of $G$, then $|G:N_{G}(H\cap N)|$ is a power of $p$.

{\bf Example 1.4.} For any group $G$, let ${\tau}(G)$ be the
set of all $CAP$-subgroups  of $G$. Then ${\tau}$  is regular
 inductive  by Lemma 1 in \cite{{Ezq1}}.

{\bf Example 1.5.} A subgroup $H$ of $G$ is called completely
$c$-permutable \cite{GS8} provided for any two subgroups  $A\leq E$
of $G$, where $H\leq E$, there is an element $x\in E$ such that
$HA^{x}=A^{x}H$. Let ${\tau}(G)$ be the set of all   completely
$c$-permutable subgroups of $G$. Then in view of \cite[Lemma
2.1(3) and Corollary 2.2(1)]{GS8}, the functor $\tau$ is hereditary inductive. Now let $H$ be a
completely $c$-permutable $p$-subgroup of $G$ and $N$ any abelian
minimal normal subgroup of $G$. For any prime $q\ne p$ dividing
$|G|$, there is a Sylow $q$-subgroup $Q$ of $G$ such that $HQ=QH$.
Then $H\cap N\ne 1$, $H\cap N=HQ \cap N$ and so $q$ does not
divide $|G:N_{G}(H\cap N)|$. Hence $\tau$ is quasiregular.

{\bf Example 1.6.} Let   ${\tau}(G)$ be  the set of all normal or
of all $S$-quasinormal subgroups of $G$, for any group $G$. Then by
\cite[1.2.7 and  1.2.8]{prod}), the funtor $\tau$  is hereditary
inductive and regular
(see \cite[1.2.14]{prod} and Example 1.4).

{\bf  Example 1.7.}  We say, by analogy with $CAP$-subgroups,
that a subgroup $H$ of $G$ is a \emph{$CAMP$-subgroup} of $G$ if for every subgroups $K\leq L\leq
G$, where $K$ is a maximal subgroup of $L$, $H$ either covers the
pair $(K, L)$ (that is, $HK=HL$) or avoids this one (that is, $H\cap
K=H\cap L$).  Let ${\tau}(G)$ be the set of all
$CAMP$-subgroups of $G$ for any group $G$. Then ${\tau}$ is
hereditary inductive by \cite[Lemma 2.3]{GuoS176}. Now let $H$ be a
$p$-subgroup of a primitive group $G$ which is a $CAMP$-subgroup of
$G$. Then $H$ is subnormal in $G$ by \cite[(Lemma 2.5)]{GuoS176}.
Let $N$ be a minimal normal subgroup of $G$. Suppose that $L=H\cap
N\ne 1$. Then $N$ is a $p$-group, so $NH$ is a subnormal
$p$-subgroup of $G$. Now let $M$ be a maximal subgroup of $G$ such
that $M_{G}=1$. Then $G=N\rtimes M$ and $H$ either covers or avoids
the pair $(M, G)$. But since $L=H\cap N\ne 1$, $H\nleq M$ and so
$G=HM$. On the other hand, $NH\cap M=1 $ by
\cite[Lemma  7.3.16]{LennoxS} since $M_{G}=1$. Therefore $NH=N=H$.
This shows that $\tau$ is $\Phi$-regular. It is not difficult to find a example
which show that $\tau$ is not quasiregular.

{\bf Example 1.8.} Recall that $H$ is said to be
\emph{$S$-quasinormally (respectively subnormally) embedded}
\cite{Bal1} in $G$ if every Sylow subgroup of $H$ is also a Sylow subgroup
of some $S$-quasinormal (respectively subnormal) subgroup of $G$.
Note that in view of Kegel's result \cite{Keg}, every
$S$-quasinormal subgroup is subnormal, so every $S$-quasinormally
embedded subgroup is also subnormally embedded.
If ${\tau}(G)$ is the set of all $S$-quasinormally embedded
subgroups of $G$ for any group $G$, then ${\tau}$ is a hereditary inductive by
\cite[Lemma 1]{Bal1} and ${\tau}$ is quasiregular (see Example 1.5
and \cite[1.2.19]{prod}). It is clear that this functor is not
regular since every Sylow subgroup is $S$-quasinormally emebedded.

{\bf Example 1.9.}   Let ${\tau}(G)$ be the set of all
modular subgroups of $G$ for any group $G$. Then ${\tau}$ is
hereditary by \cite[p.~201]{Schm}.  From Theorem 5.2.5 in \cite{Schm} it easily
follows that ${\tau}$ is   regular.

{\bf Example 1.10.} A subgroup $H$ of $G$ is called $SS$-quasinormal
in $G$ \cite{Li7} provided there is a subgroup $B$ of $G$ such that
$HB=G$ and $H$ permutes with all Sylow subgroups of $B$. If
${\tau}(G)$ is  the set of all $S$-quasinormal subgroups of $G$  for
any group $G$, then the functor $\tau $ is hereditary inductive by \cite[Lemma
2.1]{Li7} and it is also regular by \cite[Lemma 7.1(6)]{Chen-Guo}.

In what follows,  $\tau$ is always supposed to be an
inductive subgroup functor.

Now we can state  our  first result.

{\bf Theorem 1.11.}  {\sl A soluble group $G$ is supersoluble if and only if $G$
has a normal subgroup $E$ with supersoluble quotient $G/E$ such that
every maximal subgroup of every  Sylow subgroup of $F(E)$
is ${\frak{U}}_{\tau}$-supplemented
in $G $ for some $\Phi$-regular subgroup functor $\tau$. }

In view of the  remarks after Definition 1.2, Theorem 1.11  gives
the positive answer to the second part of Question A.

If $1\in \frak{F}$, then we write $G^{\frak{F}}$ to denote the
intersection of all normal subgroups $N$ of $G$ with  $G/N\in
\frak{F}$. The class $\cal F$ is said to be a \emph{formation } if
either ${\cal F}= \varnothing $ or $1\in \frak{F}$ and every
homomorphic image of $G/G^{\frak{F}}$  belongs to $ \frak{F}$ for
any group $G$. The formation $\frak{F}$ is said to be
\emph{saturated} if $G\in \frak{F}$ whenever $G/\Phi (G) \in
\frak{F}$.

In fact, Theorem 1.11  is a special case of the
next our result.

{\bf Theorem A.}  Let $\cal F$ be a saturated formation
containing $\frak{U}$ and  $X\leq E$   normal subgroups of $G$ such
that $G/E\in \frak {F}$. Suppose that every maximal subgroup of
every non-cyclic Sylow subgroup of $X$ is ${\frak{U}}_{\tau}$-supplemented
in $G$ for some $\Phi$-regular hereditary or regular subgroup functor $\tau$ such that every
$\tau $-subgroup of $G$  contained in $X$ is subnormally embedded in $G$.
If   $X=E$ or $X=F^{*}(E)$, then  $G\in \frak{F}$. Moreover, in the case when $\tau$ is regular, then $E\leq
Z_{\frak{U}}(G)$.

The following theorem is an analogue  of the previous one. But the
methods of their proofs are quite different (see
Sections 3 and 4).

{\bf  Theorem B.}   {\sl Let $\frak{F}$ be a  saturated formation
containing $\frak{U}$ and $X\leq  E$ a normal subgroup of $G$ such
that $G/E \in \frak{F}$. Suppose that  for every non-cyclic Sylow
subgroup $P$ of $X$  every cyclic subgroup   of $P$ of
 prime order or  order 4 (if $P$ is a non-abelian group) is  ${\frak{U}}_{\tau}$-supplemented
 in $G $. Suppose that at least one of the following holds:}

(i) {\sl $\tau$ is hereditary $\Phi$-quasiregular  and $X=E$};

(ii) {\sl $\tau $ is hereditary quasiregular, $E$ is soluble
and $X=F(E)$};

(iii) {\sl $\tau $ is regular, and $X=F^{*}$ or $X=E$}.

{\sl Then  $G\in \frak{F}$.}

By analogy with Theorem 1.11,  Theorems A and B cover and unify
the results in many papers. Some of them we discuss in   Section 5.

Finally, note that Theorem A and B remain to be new for each
concrete subgroup functor $\tau $, for example, if we supposed that
$\tau $ is one of the functors in Examples 1.4--1.10.

All unexplained notation and terminology are standard. The reader is
referred to \cite{100},  \cite{DH},    \cite{Guo}, or  \cite{BalE}, if
necessary.

\section{Base  lemmas}

{\bf Lemma 2.1.} {\sl   Let $ \frak{F}$ be  a
saturated formation,  $K\leq H\leq G$ and $N\trianglelefteq G$.
Suppose that the pair $(K, H)$ satisfies the  $\frak{F}$-supplement
condition in $G$. }

(1) {\sl If either $N\leq H$  or $(|H|, |N|)=1$, then the pair
$(KN/N, HN/N)$ satisfies the  $\frak{F}$-supplement condition in
 $G/N$.}

(2) {\sl If $H\leq E\leq  G$ and $\frak{F}$ is  hereditary, then the
pair $(K, H)$ satisfies the ${\cal F }$-supplement condition  in
$E$.}

(3) {\sl If  $K\leq V\leq H$, then  the pair $(V, H)$  satisfies the
 $\frak{F}$-supplement condition   in $G$.}

{\bf Proof.} Let $T$ be a   subgroup of $G$ such that
 $HT=G$ and $H\cap T\subseteq  KZ_\frak{F}(T)$.

(1) Clearly, $(HN/N)(TN/N)=G/N$.  Moreover, $HN\cap HN=(H\cap T)N$.
Indeed, if either $N\leq H$ or $N\leq T$, it is clear. But if
$N\nleq H$, we have $N\leq T$ since in this case $(|H|, |N|)=1$ by
hypothesis. Hence $$(HN/N)\cap (NT/N)=(NH\cap NT)/N= N(H\cap
T)/N\subseteq N(KZ_\frak{F}(T))/N= (NK/N)(Z_\frak{F}(T)N/N)$$ $$ \leq
(KN/N)(Z_\frak{F}(TN/N))$$ since $Z_\frak{F}(T)N/N)\leq
Z_\frak{F}(TN/N)$ by \cite[Lemma 2.2(4)]{GuoS-JA}.  Hence the pair $
(KN/N, HN/N)$ satisfies the $\frak{F}$-supplement condition
 in $G/N$.

(2)  Let $T_{0}=T\cap E$. Then  $E=E\cap HT=H(T\cap E)= HT_{0}$.
Moreover,   $$H \cap T_{0} =H\cap T\subseteq  KZ_\frak{F}(T)\leq
K(Z_\frak{F}(T)\cap E)$$ $$\leq KZ_{\cal F}(T\cap E)= KZ_{\cal
F}(T_{0})$$  by \cite[Lemma 2.2(5)]{GuoS-JA}. This shows that the pair
$(K, H)$ satisfies the  $\frak{F}$-supplement condition  in $E$.

(3) This is clear.

{\bf  Lemma 2.2.} {\sl  Let $\frak{F}$ be a   saturated formation.
Let  $H\leq G$ and  $N$ be  a normal subgroup of $G$.

(1) If  $N\leq H$, then  $H/N$ is  $\frak{F}_{\tau
}$-supplemented in $G/N$ if and only if  $H$ is   ${\cal F}_{\tau
}$-supplemented  in $G$.

(2)  Suppose that $S$ is a $\tau$-subgroup of $G$ such that
$H_G\leq S\leq H$  and the pair $(S, H)$  satisfies the
$\frak{F}$-supplement condition in $G$. Then $H$ is $\frak{F}_{\tau}$-supplemented  in $G$.

(3) If $H$ is $\frak{F}_{\tau}$-supplemented in $G$, and either $N\leq H$ or $(|H|, |N|)=1$, then $HN/N$ is
$\frak{F}_{\tau }$-supplemented  in $G/N$.

(4) If $H$ is $\frak{F}_{\tau}$-supplemented in $G$, $\frak{F}$ is  hereditary,  $H\leq E\leq G$ and  $\tau$
is hereditary, then $H$ is $\frak{F}_{\tau}$-supplemented in $E$. }

{\bf Proof.}  For any subgroup $V$   of $G$ we  put
 $\bar{V}=VN/N$  and $\widehat{V}=VH_{G}/H_{G}$.
 Let $f$  be the  canonical isomorphism
from $(G/N)/(H_{G}/N)$ onto $G/H_{G}$.

Let $\bar{T}/{\bar{H}}_{\bar{G}}$ be a  subgroup of
$\bar{G}/{\bar{H}}_{\bar{G}}$ such that
$\bar{G}/\bar{H}_{\bar{G}}=(\bar{H}/\bar{H}_{\bar{G}})(\bar{T}/\bar{H}_{\bar{G}})$
and
$$(\bar{H}/{\bar{H}}_{\bar{G}}) \cap ( \bar{T}/{\bar{H}}_{\bar{G}})
 ) \subseteq (\bar{S}/{\bar{H}}_{\bar{G}}) Z_{\cal
F}(\bar{T}/{\bar{H}}_{\bar{G}})$$ for some $\tau$-subgroup
$\bar{S}/{\bar{H}}_{\bar{G}}$ of $\bar{G}/{\bar{H}}_{\bar{G}}$
contained in $\bar{H}/{\bar{H}}_{\bar{G}}$.

Then
$$f((\bar{H}/{\bar{H}}_{\bar{G}}) \cap ( \bar{T}/{\bar{H}}_{\bar{G}})
)  \leq
f(\bar{S}/{\bar{H}}_{\bar{G}})f(Z_{\cal
F}(\bar{T}/{\bar{H}}_{\bar{G}}).$$

Note that $$f((\bar{H}/{\bar{H}}_{\bar{G}}) \cap (
\bar{T}/{\bar{H}}_{\bar{G}}) ) )=\widehat{ H}\cap \widehat {T}$$ and
${\bar{H}}_{\bar{G}}=(H/N)_{G/N}=H_{G}/N$.

On the other hand,
$$f(\bar{S}/{\bar{H}}_{\bar{G}})=S/H_{G}=\widehat{S},$$
where $\widehat{S}$ is a $\tau$-subgroup of $\widehat{G}$ since
 $\bar{S}/{\bar{H}}_{\bar{G}}$ is a $\tau$-subgroup of
$\bar{G}/{\bar{H}}_{\bar{G}}$.

Similarly, $f(Z_\frak{F}(\bar{T}/{\bar{H}}_{\bar{G}}))
=Z_\frak{F}(\widehat{ T})$. Therefore $H$ is $\frak{F}_{\tau
}$-supplemented in $G$.

Now suppose that $H$ is  $\frak{F}_{\tau }$-supplemented  in $G$. Then
by considering the canonical isomorphism $f^{-1}$ from $G/H_{G}$
onto $(G/N)(H_{G}/N)$, one can prove analogously that $H/N$ is
$\frak{F}_{\tau }$-supplemented  in $G/N$. The second assertion of
(1) can be proved similarly.

(2) By the hyperthesis, the pair $(S, H)$ satisfies the
$\frak{F}$-supplement condition  in $G$. Hence by Lemma 2.1(1),
$(S/H_{G}, H/H_{G})$ satisfies the  $\frak{F}$-supplement condition
in $G/H_G$. Thus  $H$ is $\frak{F}_{\tau }$-supplemented  in $G$.

(3)  Let $\widehat{S}$ be a $\tau$-subgroup of $\widehat{G}$
contained in $\widehat{H}$ such that the pair $(\widehat{S},
\widehat{H}) $ satisfies the  $\frak{F}$-supplement condition
 in $\widehat{G}$.

Then the pair $(\widehat{S}\widehat{N}/\widehat{N},
\widehat{H}\widehat{N}/\widehat{N}) $ satisfies the    ${\cal
F}$-supplement condition  in $\widehat{G}/\widehat{N}$  by Lemma 2.1(1).
Let $h$  be the  canonical isomorphism from
$(G/H_{G})/(H_{G}N/H_{G})$ onto $G/NH_{G}$. Then
$h(\widehat{S}\widehat{N}/\widehat{N})=SN/NH_{G}$ and
$h(\widehat{H}\widehat{N}/\widehat{N})=HN/NH_{G}$. Hence the pair
$(SN/NH_{G}, HN/NH_{G} )$ satisfies the  $\frak{F}$-supplement
condition  in $G/NH_{G}$. Note also that $SN/NH_{G} $ is a
$\tau$-subgroup of $G/NH_{G}$ since $\tau$ is inductive and
$\widehat{S}$ is a $\tau$-subgroup of $\widehat{G}$. Hence,
$HN/NH_{G} $ is  $\frak{F}_{\tau}$-supplemented in $G/NH_{G}$ and so
$(HN/N)/(H_{G}N)/N$ is  $\frak{F}_{\tau }$-supplemented in
$(G/N)/(NH_{G}/N)$, which implies that $HN/N$ is $\frak{F}_{\tau
}$-supplemented   in $G/N$ by (1).

(4) By hypothesis, for some $\widehat{S}\leq \widehat{H}$, where
 $\widehat{S}$ is $\tau$-subgroup of $\widehat{G}$,  the pair $(\widehat{S} ,
 \widehat{H})$ satisfies the  $\frak{F}$-supplement condition  in $
 \widehat{G}$. Then, by Lemma 2.1(2), the pair $(\widehat{S} ,
 \widehat{H})$ satisfies the  $\frak{F}$-supplement condition  in $
 \widehat{E}$. Hence, by Lemma 2.1(3), the pair $(\widehat{S}\widehat{H_{E}}/\widehat{H_{E}},
 \widehat{H}/\widehat{H_{E}} )$
satisfies the  $\frak{F}$-supplement condition in $
 \widehat{E}/\widehat{H_{E}}$, where $(\widehat{S}\widehat{H_{E}}/\widehat{H_{E}}$  is $\tau$-subgroup of $\widehat{E}/\widehat{H_{E}}$.
Hence $SH_{E}/H_{E}$ is a $\tau$-subgroup of $E/H_{E}$ and the pair
$(SH_{E}/H_{E}, H/H_{E})$ satisfies the  $\frak{F}$-supplement
condition  in $E/H_{E}$. This shows that $H$ is $\frak{F}_{\tau}$-supplemented in
$E$.

\section{Proof of Theorem A}

The following lemma  is a corollary of \cite[IV,  (6.7)]{DH}.

{\bf  Lemma 3.1}  {\sl Let $\frak{F}$ be a  saturated formation and
$F$  the canonical local  satellite of $\frak{F}$} (See \cite[Lemma
2.7]{Guo181} or \cite[p. 361]{DH}).  {\sl  Let $P$ be a normal $p$-subgroup
of $G$. If $P/\Phi (P)\leq Z_{\frak{F}}(G/\Phi (P))$, then $P \leq
Z_{\frak{F}}(G)$. }

{\bf Lemma 3.2} (See Lemma  2.14 in \cite{Guo181}).  {\sl  Let
${\cal F}$ be a saturated  formation and  $F$ the canonical local
 satellite of $\cal F$.
 Let $E$ be a
normal $p$-subgroup of  $G$. Then   $ E\leq Z_{\cal F}(G)$ if and
only if  $G/C_{G}(E)\in F(p)$.}

{\bf  Lemma 3.3.} {\sl  Let $\frak{F}$ be a  saturated formation and
$G=PT$, where $P$ is a normal $p$-subgroup of $G$.  If $P\cap
Z_{\frak{F}}(T)$ is normal in $P$, then $P\cap Z_{\frak{F}}(T)\leq
Z_{\frak{F}}(G)$. }

{\bf Proof.} First note that since   $G=PT$   and $P\cap
Z_{\frak{F}}(T)$ is normal in $P$, $P\cap Z_{\frak{F}}(T)$ is normal
in $G$.  Let $F$ be the canonical local satellite of
$\frak{F}$ and   $H/K$ a chief factor of $G$ below $P\cap
Z_{\frak{F}}(T)$. Then $T/C_{T}(H/K)\in F(p)$ by Lemma 3.2. Hence
$G/C_{G}(H/K)=PT/C_{G}(H/K)=PT/P(C_{G}(H/K)\cap T)\simeq
T/P(C_{G}(H/K)\cap T)\cap T=T/(C_{G}(H/K)\cap T)=T/C_{T}(H/K)\in
F(p)$ by \cite[Lemma 2.11]{SkJGT}  or   \cite[A, (10.6)(b)]{DH}. Hence
$P\cap Z_{\frak{F}}(T)\leq Z_{\frak{F}}(G)$.

{\bf  Lemma 3.4.} {\sl Let $L\leq  V\trianglelefteq P$, where $P$ is
a Sylow $p$-subgroup of a group $G$, and  $N$ and $M$ are different
normal subgroups of $G$.  Suppose that $|G/M: N_{G/M}(( LM/M)\cap
(NM/M))$ is a power of  $p$. Then:}

(1) {\sl $(( LM/M)\cap (NM/M))^{G/M}\leq VM/M$.}

(2) {\sl If $N$ is a non-abelian  minimal normal subgroup of $G$, then $L\cap N=1$.}

(3) {\sl  If $NL\cap M=1$, then $(L\cap N)^{G}\leq VM$}

{\bf Proof.} (1) It is clear that $LM/M \leq  VM/M \trianglelefteq
PM/M$ where $PM/M$ is a Sylow $p$-subgroup of $G/M$. On the other
hand, since   $|G/M: N_{G/M}(( LM/M)\cap (NM/M))$ is a power of $p$,
we have $(( LM/M)\cap (NM/M))^{G/M}=(( LM/M)\cap (NM/M))^{N_{G/M}(((
LM/M)\cap (NM/M)))(PM/M)}=(( LM/M)\cap (NM/M))^{PM/M}$. Hence we
have (1).

(2) Suppose that  $B=L\cap N\ne 1.$ Then $( LM/M)\cap (NM/M)\ne 1$
and so $N\simeq NM/M \leq  (( LM/M)\cap (NM/M))^{G/M}\leq  PM/M$,
which implies that $N$ is a $p$-group.

(3) Since $NL\cap M=1$, $(LM/M)\cap (NM/M)= (L\cap N)M/M$. On the
other hand, $((L\cap N)M/M)^{G/M}=(L\cap N)^{G}M/M$. Hence
(3) is a corollary of (1).

A normal subgroup $N$ of $G$ is said to be
\emph{${\frak{F}}\Phi$-hypercentral}  in $G$ \cite{shem-skJA2009} if
either $N=1$ or $N\ne 1$ and every non-Frattini  chief factor of $G$
below $N$ is ${\cal F}$-central in $G$. The product of all normal
${\frak{F}}\Phi$-hypercentral subgroups is   denoted by
$Z_{\frak{F}\Phi}(G)$ \cite{shem-skJA2009}.

{\bf  Proposition 3.5.}  {\sl Let $\frak{F}$ be a  saturated
formation containing all supersoluble groups, $\tau$ be $\Phi$-quasiregular (quasiregular, respectively) and  $P$  a
non-identity normal $p$-subgroup of  $G$. Suppose that every maximal
subgroup of $P$ is ${\frak{F}}_{\tau}$-supplemented in $G $. Then
$P\leq Z_{\frak{F}\Phi}(G)$ ($E\leq Z_{\frak{F}}(G)$, respectively).}

{\bf Proof.}  Suppose that the proposition is false and let $(G,
P)$ be a counterexample with $|G|+|P|$ minimal.  Let
$Z=Z_{\frak{F}\Phi}(G)$ ($Z=Z_{\frak{F}}(G)$, respectively). Let
$G_{p}$ be a Sylow $p$-subgroup of $G$.

(1) {\sl   $P\not \leq Z_{\frak{U}\Phi}(G)$ ($P\not \leq
Z_{\frak{U}}(G)$, respectively). } (This follows from the hypothesis
that $\frak{F}$  contains all supersoluble groups and the choice of
$G$ ).

(2) {\sl $P$ is not a minimal normal subgroup of $G$.}

Suppose that $P$ is a minimal normal subgroup of $G$. Then $P\cap
Z=1$. Let $H$ be a maximal subgroup of $P$ such that $H$ is normal
in $G_{p}$. Then $H_{G}=1$. Let $S\in \tau (G)$ and $T$ be subgroups
of $G$ such that $S\leq H$, $HT=G$ and $H\cap T\subseteq
SZ_{\frak{F}}(T)$. Suppose that $T\ne G$. Then $1 < P\cap T < P$,
where $P\cap T$ is normal in $G$ since $P$ is abelian, which
contradicts the minimality of $P$. Hence $T=G$, so $H=H\cap  T\leq
SZ$, which implies that $H=S(H\cap Z)=S$ is a $\tau$-subgroup of
$G$. It is clear that $P\nleq \Phi (G)$. Hence for some maximal
subgroup $M$ of $G$ we have $G=P\rtimes M$. Since  $\tau$ is
$\Phi$-quasiregular and $HM_{G}/M_{G}$ is normal in
$G_{p}M_{G}/M_{G}$, $HM_{G}/M_{G}$ is normal in $G/M_{G}$ by Lemma 3.4.
Hence $HM_{G}/M_{G}=PM_{G}/M_{G}$, which implies that $H=P$, a
contradiction. Hence we have (3).

(3)  {\sl If $N$ is a minimal normal subgroup of $G$ contained in
$P$, then $P/N\leq Z_{\frak{F}\Phi}(G/N)$ ($P/N\leq
Z_{\frak{F}}(G/N)$, respectively) and $Z\cap P=1$.}

Indeed, by Lemma 2.2      the hypothesis holds for $(G/N, P/N) $. Hence
$P/N \leq Z_{\frak{F}\Phi}(G/N)$ ($P/N \leq Z_{\frak{F}}(G/N)$,
respectively) by the choice of $(G, E)$. Hence $N\nleq Z$ by
\cite[Lemma 2.2]{shem-skJA2009} and the choice of $(G, P)$.

(4) {\sl $P\leq Z_{\frak{F}\Phi}(G)$.}

Suppose that $P\nleq Z_{\frak{F}\Phi}(G)$. Then, in view of (3) and
\cite[Lemma 2.2]{shem-skJA2009}, $\Phi (G)\cap P=1$. Hence  $P=N\times D$
for some normal subgroup $D$ of $G$ by \cite[Lemma 2.15]{Guo181}, where
$D\ne 1$ by (2). Let $R$ be a minimal normal subgroup  of $G$
contained in $D$. Then, by \cite[A, 9.11]{DH},  $RN/N\nleq \Phi
(G/N)$. Hence in view of (3) and the $G$-isomorphism $R\simeq RN/N$,
$R$ is $\frak{F}$-central in $G$, and so $P\leq Z_{\frak{F}\Phi}(G)$ by
(3) and \cite[Lemma 2.2]{shem-skJA2009}. Hence we have (4).

(5) {\sl $\tau$ is quasiregular} (This follows from (4) and the
choice of $(G, P)$).

(6)  {\sl If $N$ is a minimal normal subgroup of $G$ contained in
$P$, then $N$ is the unique minimal normal subgroup of $G$ contained
in $P$} (see the proof of (3)).

(7)  $\Phi (P)\ne 1$.

Suppose that $\Phi (P)= 1$. Then $P$ is an elementary abelian
$p$-group. Let $W$ be a maximal  subgroup of $N$ such that  $W$ is
normal  in   $G_p$. We show that $W$ is normal in $G$. Let $B$
be a complement of $N$ in $P$ and   $V=WB $. Then  $V$ is a maximal
subgroup of $P$ and, evidently, $V_{G}=1$ by (6).

Let $S\in \tau (G)$ and $T$ be subgroups of $G$ such that $S\leq V$,
$VT=G$ and $V\cap T\subseteq SZ_{\frak{F}}(T)$. Assume that $T=G$. Then
$V=V\cap T\leq SZ$ and so $V= S(V\cap Z)$. But in view of (3),
$Z\cap P=1$. Hence $V= S$ and thereby $W= WB\cap N=V\cap N=S\cap N$.
 Since ${\tau}$ is quasiregular by (5) and $W$ is
normal in $G_{p}$, we have that that  $W$ is normal in $G$ by Lemma 3.4.
Let $T\ne G$. Then $1\ne T \cap P < P$. Since $G=VT=PT$ and $P$ is
abelian, $T\cap P$ is normal in $G$. Hence $N\leq T$, which implies
that $W\leq T\cap V\subseteq SZ_\frak{F}(T)\cap P=S(Z_\frak{F}(T)\cap
P)$. Since $P$ is abelian, $Z_\frak{F}(T)\cap P$ is normal in $P$,
so $Z_\frak{F}(T)\cap P\leq Z\cap P=1$ by Lemma 3.3. This implies that
 $W= S\cap N$ and so $W$ is normal in $G$ by Lemma 3.4.

Finally,  as above, we get the same conclusion in the case when
$N\leq T$. Therefore $W$ is normal in $G$ and so $W=1$, which
implies that  $|N|=p$. This contradiction shows that we have $\Phi
(P)\ne 1$.

{\sl The final contradiction.}

By (7), $\Phi (P)\ne 1$. Let $N$ be a minimal normal subgroup of $G$
contained in $\Phi (P).$ Then $P/N\leq Z_{\frak{F}}(G/N)$ by (3). It
follows that $P/\Phi (P)\leq Z_{\frak{F}}(G/\Phi (P))$. Thus $P\leq
Z$ by Lemma 3.1. This contradiction completes the proof.

{\bf  Theorem 3.6.}   {\sl Let $\frak{F}$ be the class of all
$p$-supersoluble groups. Let   $E$ be  a normal subgroup of  $G$ and
$P$ a Sylow  $p$-subgroup of $E$ of order $|P|=p^{n}$, where $n > 1$
and $(|E|, p-1)=1$.  Suppose that $\tau $ is  $\Phi$-regular and
every $\tau$-subgroup  of $G$ contained in $P$  is subnormally
embedded in $G$. If  every maximal subgroup of $P$ is
${\frak{F}}_{\tau}$-supplemented in $G $, then $E$ is $p$-nilpotent.}

{\bf Proof.}   Suppose that this theorem is false and let $(G, E)$
be a counterexample with $|G|+|E|$ minimal. We proceed via the following steps.

(1) $O_{p'}(E)=1$.

Suppose that $O_{p'}(E)\ne 1$. Since $O_{p'}(E)$ is characteristic
in $E$, it is normal in $G$ and the hypothesis holds for
$(G/O_{p'}(E), E/O_{p'}(E))$ by Lemma 2.2. The choice of $G$ implies
that $E/O_{p'}(E)$ is $p$-nilpotent. It follows that $E$ is
$p$-nilpotent, a contradiction.

(2) {\sl If $O_{p}(E)\ne 1$, then $E$ is $p$-soluble.}

Indeed, by Lemma 2.2, the hypothesis holds for $(G/O_{p}(E),
E/O_{p}(E))$. Hence in the case when $O_{p}(E)\ne 1$, $E/O_{p}(E)$
is $p$-nilpotent by the choice of $(G, E)$, which implies the
$p$-solubility of $E$.

(3) {\sl $O_{p}(E)\ne 1$.}

Suppose that  this is false. Then, in view of
\cite[Lemma 3.4(3)]{Guo181}, for any subnormal subgroup $L$ of $G$
contained in $E$ we have neither $L$ is a p-group  nor  a
$p'$-group. Let  $N$ be a minimal normal subgroup of $G$ contained
in $E$. Then  $N$ is non-abelian group. Then since $(|E|, p-1)=1,$ we have that $p=2$ by Feit-Thompson theorem. It is clear that $|N_2|>2$.

We claim that for any minimal normal subgroup $L$ of $G$ contained in $E$
and  for any $\tau$-subgroup $S$ of $G$ we have $S\cap L=1$. Indeed, assume that $S\cap L\ne 1$. Let $M
$ be a maximal subgroup of $G$ such that $LM=G$.
Since ${\tau}$ is $\Phi$-regular,  $|G/M_{G}: N_{G/M_{G}}(( SM_{G}/M_{G})\cap
(LM_{G}/M_{G}))|$ is a power of  $2$. Then  $L$ is abelian by
Lemma 3.4(2) and so $L\leq O_2(E)=1,$ a contradiction.

Let $H$ be an arbitrary maximal subgroup of $P$. It is clear that
$H_{G}=1$. Hence   there exists a subgroup $T$ such that $G=HT$ and
$H\cap T\subseteq SZ_{\cal F}(T)$ for some ${\tau }$-subgroup $S$ of $G$
contained in $H$.

Suppose that $S\ne 1$. Let  $W$ be a subnormal subgroup of $G$ such
that  $S$ is a Sylow $2$-subgroup of $W$. In view of
\cite[Lemma 3.4(2)]{Guo181}, we may assume, without loss of generality,
that $W\leq E$.  Let $L$ be a minimal subnormal
subgroup of $G$ contained in $W$. Then  $L$ neither  is a 2-group
nor a $2'$-group. Therefore $L$ is a non-abelian simple group and
$L_{2}=S\cap L$ is a Sylow $2$-subgroup of $L$ since $S$ is a Sylow
$2$-subgroup of $W$. It is clear that  $R=L^{G}$ is a minimal normal
subgroup of $G$ and $S\cap R\ne 1$, contrary to (a).  Therefore $S=1$. Hence  every maximal
subgroup $H$ of $P$ has a supplement  $T$ in $G$ such that $H\cap
T\leq Z_{\cal F}(T)$.

We now show that $V=T\cap E$ is 2-nilpotent. Let $V_{2}$ be a Sylow
$2$-subgroup of $V$ containing $M\cap V$. Then $|V_{2}:V\cap M|\leq
|P:M|=2$. Therefore for a Sylow $2$-subgroup $Q$  of $VZ_{\cal
F}(T)/Z_{\cal F}(T)$
 we have $|Q| $ divides $2$. This induces that  $VZ_{\cal F}(T)/Z_{\cal F}(T)\simeq
V/V\cap Z_{\cal F}(T)$ is 2-nilpotent. It is well-known that  the
class of all $2$-nilpotent groups is a hereditary saturated
formation.  Hence in view of \cite[Lemma 2.2(5)]{GuoS-JA}, $V\cap
Z_{\cal F}(T)\leq Z_{\cal F}(V)$. Thus
 $V=T\cap E$ is
$2$-nilpotent. But $E=E\cap TM=M(T\cap E)$, so every maximal
subgroup of $P$ has a $2$-nilpotent supplement $T$ in $E$. It is
clear that a Hall $2'$-subgroup of $T\cap E$ is also a Hall
$2'$-subgroup of $E$. Therefore $E$ is $2$-nilpotent  by
\cite[Lemmas 3.7 and  3.8]{Guo181}. It follows that $N$ is a 2-group, a
contradiction. Hence we have (3)

(4) {\sl There is a maximal subgroup  $D$ of $G$ such that $ND=G$,
$D_{G}\cap E=1$  and $E=N\rtimes M$,
 where   $M=D\cap E$  and  $N=O_{p}(E)=C_{E}(N)$ is a minimal normal subgroup
of $G$ and $M$ is $p$-nilpotent.  In particular, $E$ is
$p$-soluble.}

In view of (3), $O_{p}(E)\ne 1$. Let $N$ be a minimal normal
subgroup of $G$ contained in $O_{p}(E)$.  Then the hypothesis holds
for $(G/N, E/N)$ by Lemma 2.2.
Therefore $E/N$ is
$p$-nilpotent by the choice of $(G, E)$, and so $E$ is $p$-soluble. It
follows  that $N$ is the unique minimal normal subgroup of $G$
contained in $E$. If $N\leq \Phi (G)$, then $E$  is $p$-nilpotent by
\cite[Corollary 1.6]{Guo181}. Hence $N\nleq \Phi (G)$ and so $G=N\rtimes D$
for some maximal subgroup $D$ of $G$.   Since $O_{p}(G)\leq
C_{G}(N)$ by  \cite[Lemma 2.11]{SkJGT}  or   \cite[A, 10.6(b)]{DH}, we have that
$O_{p}(G)\cap D$ is normal in $G$. Hence $O_{p}(G)\cap D\cap E$ is
normal in $G$. Note that $E=N\rtimes (D\cap E)$, so $D_{G}\cap M=1$
and
$$O_{p}(E)=O_{p}(G)\cap E= N(O_{p}(E)\cap D\cap E),$$ where
$O_{p}(E)\cap
D\cap E=O_{p}(G)\cap D\cap E$ is normal in $G$. Hence $O_{p}(E)\cap D\cap E=1$,
and so $N=O_{p}(E)$. Finally, since $E$ is $p$-soluble and $O_{p'}(E)=1$ by (1),
we have $C_{E}(N)=N$ by \cite[Chapter 6, 3.2]{Gor}.

(5) {\sl If $H/K$ is a chief factor  of $E$ below  $N$, then $|H/K|
> p$ }

By Proposition  4.13(c) in \cite[A]{DH}, $N=N_{1}\times \ldots
\times N_{t}$, where $N_{1}, \cdots , N_{t}$ are minimal normal
subgroups of $E$, and
from the proof of this proposition we see  that $|N_{i}=|N_{j}|$
for all $i, j \in \{1, \cdots ,t\}$. Hence for any
chief factor $H/K$  of $E$ below  $N$ we have $|H/K|=|N_{1}|$ by
\cite[Lemma 2.2]{shem-skJA2009}. Suppose that $|H/K|=p$.
 Since $(p-1, |E|)=p$,
$C_{E}(H/K)=E$. Hence $N\leq Z_{\infty }(E)$, which implies the
$p$-nilpotency of $E$ by (4). This contradiction shows that (5)
holds.

(6) {\sl If a non-identity subgroup $S$ of $P$ is subnormally
embedded in $G$, then $S\cap N\ne 1$.}

Indeed, let $ W$ be a subnormal subgroup of  $G$  such that  $S$ is
a Sylow   $p$-subgroup of $W$. If  $S\cap N=1$, then $W\cap N=1$.
Hence by (4) and \cite[Lemma 3.4(4)]{Guo181},
$W\leq C_{E}(N)=N$, a contradiction. Thus  (6) holds.

(7)  {\sl $M=N_{E}(M_{p'})$,  where  $M_{p'}$ is the Hall
$p'$-subgroup of $M$.}

Let $J=N_{E}(M_{p'})$.  Suppose that $M < J$. Then $J=J\cap
NM=M(J\cap N)$ and therefore $J\cap N\ne 1$. Since $E=NJ$, $J\cap N$
is normal in $E$ and $E/C_{E}(J\cap N)$ is a $p$-group. If $F$ is
the canonical local satellite of the saturated formation of all
nilpotent, the $F(p)$ is the class of all $p$-groups by (see
\cite[IV]{DH}). Hence in view of Lemma 3.2, $J\cap N\leq Z_{\infty}(E)$.
Hence for some minimal normal subgroup
$C$ of $E$ contained in $N$ we have $|C|=p$, which contradicts (5).

{\sl Final contradiction.}

Let $M_{p}\leq D_{p}$, where $M_{p}$ is a Sylow $p$-subgroup of $M$
and $D_{p}$ is a Sylow $p$-subgroup of $D$.  Without loss of
generality, we may suppose that  $M_{p}\leq P$. Then  $NM_{p}= P$
and $ND_{p}$  is a Sylow $p$-subgroup of $G$.
Let $N_0\leq N$ be a normal subgroup of   $ND_{p}$ such
that $|N:N_0|=p$.  Let $W=N_0D_{p}$  and $V=N_0M_{p}$. Then
$W$ is maximal in $ND_{p}$ and
$V$ is maximal in $P$ such that $V_{G}=1$.

(i) {\sl For any $\tau $-subgroup $S$ of $G$ contained in $V$, we
have $S\cap N=1$.}

Assume that this is false. It is cleat that $SN\cap D_{G}=1$. Since ${\tau}$ is $\Phi$-regular, $N\leq (S\cap N)^{G}D_{G}\leq
WD_{G}$ by Lemma 3.4, Hence $ N=N\cap N_0D_{p}D_{G}=N_0(N\cap
D_{p}D_{G})=N_0$. This contradiction shows that we have (i).

(ii) {$V$ has no a $p$-nilpotent  supplement in $E$.}

Assume that $V$ has a $p$-nilpotent  supplement $T_{0}$ in $E$. Then
a Hall $p'$-subgroup $T_{p'}$ of $T_{0}$ is a  Hall $p'$-subgroup
of $E$. By (4), $E$ is $p$-soluble and so any two   Hall $p'$-subgroup
of $E$ are conjugate in $E$. Without loss of generality, we may  assume that $T_{p'}\leq M$, so  $T_0\leq
M$ by (7). It follows that $E=VT_{0}=VM$. But since $M_{p}\leq V$
and $V$ is maximal in $P$, $VM\ne E$. This  contradiction shows that
we have (ii).

By hypothesis, $V$ is ${\frak{F}}_{\tau}$-supplemented in $G$, so there exists a
subgroup $T$ such that $G=VT$ and $V\cap T\subseteq SZ_{\cal F}(T)$ for
some ${\tau }$-subgroup $S$ of $G$ contained in $V$.

Assume that $S\ne 1$.  Then $S\cap N\ne 1$ by (6), contrary to (i).
Hence  $S=1$, so
$$
V\cap T_{0}=V\cap T \leq Z_{\cal F}(T) \cap T_{0} \leq Z_{\cal
F}(T_{0}) $$ by \cite[Lemma 2.2(5)]{GuoS-JA}.  Hence, as in the proof of
(3), one can show that $T_{0}=T\cap E$ is  $p$-nilpotent. But this
contradicts (ii).

The theorem is proved.

{\bf  Corollary 3.7.}  {\sl  Let $E$ be  a non-identity normal
subgroup $G$. Suppose that every maximal subgroup of every
non-cyclic Sylow subgroup of $E$ is ${\frak{U}}_{\tau}$-supplemented
in $G $ for some regular subgroup functor $\tau$ such that every  $\tau $-subgroup of $G$  contained in $E$ is subnormally embedded in $G$.
Then $E\leq Z_{\frak{F}}(G)$.}

{\bf Proof.}  Suppose that this corollary  is false and let $(G, E)$
be a counterexample with $|G|+|E|$ minimal. Let
$p$ be the smallest prime dividing $|E|$. Then $E$ is $p$-nilpotent by Theorem 3.6 and
\cite[Chapter 7, 6.1]{Gor}.  Let $V$ be the Hall
$p'$-subgroup of $E$. Then $V$ is characteristic in $E$ and so it is
normal in $G$. If $V=1$, then $E\leq Z_{\frak{F}}(G)$ by Proposition 3.5.
Hence $V\ne 1$, and so the hypothesis holds for  $(G, V)$. The  choice of
$(G, E)$ implies that $V\leq Z_{\frak{F}}(G)$. On the other hand,
the hypothesis holds for $(G/V, E/V)$ by Lemma 2.2, so $E/V\leq
Z_{\frak{F}}(G/V)$. Therefore $E\leq Z_{\frak{F}}(G)$ by \cite[Lemma
2.2]{shem-skJA2009}. This contradiction completes the proof.

{\bf  Lemma 3.8.} {\sl Let $P$ be a normal non-identity $p$-subgroup
of $G$ with   $|P| > p$  and $P\cap \Phi (G)=1$. Suppose that $\tau$ is $\Phi$-quasiregular and every maximal subgroup of $P$
is ${\frak{U}}_{\tau}$-supplemented in $G$.
Then some maximal subgroup of $P$ is normal in
$G$.}

{\bf Proof.} Let $G_{p}$ be a Sylow $p$-subgroup of $G$ and
$Z=Z_{\frak{U}}(G)$. Since $P\cap \Phi (G)=1$, $P=N_{1}\times \cdots
\times N_{t}$, where $N_{i}$ is a minimal normal subgroup of $G$,
for all $i= 1,  \ldots t$ by  \cite[Lemma 2.15]{Guo181}.  Hence
 $|N_{i}|\ne p$ for all $i=1,  \ldots , t$. Then $P\cap
Z=1$ and $t> 1$ (see (2) in the proof of Proposition 3.5). Moreover, the hypothesis
holds for $(G/N_{1}, P/N_{1})$ by Lemma 2.2 and \cite[Lemma 2.2]{shem-skJA2009},
 so by induction some maximal
subgroup $M/N_{1}$ of $P/N_{1}$ is normal in $G/N_{1}$. Then a maximal
subgroup $M$ of $P$ is normal in $G$.

{\bf  Lemma 3.9}  {\sl Let $E$ be a normal subgroup of $G$ and $\tau
$  a $\Phi$-regular inductive subgroup functor such that every
primary $\tau$-subgroup of $G$ is subnormally embedded in $G$. If
every  maximal subgroup of every non-cyclic Sylow subgroup of $E$ is
${\frak{U}}_{\tau}$-supplemented in $G $, then $E$ is
supersoluble.}

{\bf Proof. } Suppose that this lemma is false and let $(G, E)$ be a counterexample  with $|G|+|E|$ minimal.
Let $P$ be a Sylow $p$-subgroup of $E$ where $p$ is
the smallest prime dividing $|E|$. By Theorem 3.6, $E$ is $p$-nilpotent.
Let $V$ be the Hall $p'$-subgroup of $E$. Then $V$ is normal in $G$ and the  hypothesis
holds for $(G, V)$. Hence $V$ is supersoluble  by the choice of $(G,
E)$. Then a Sylow $q$-subgroup $Q$  of $V$, where $q$ is the
largest prime dividing $|V|$, is normal and so it is characteristic in
$V$. Hence $Q$ is normal in $G$ and the hypothesis holds for $(G/Q,
E/Q)$ by Lemma 2.2. The choice of $(G, E)$  implies that $E/N$ is
supersoluble. On the other hand, by Proposition 3.5, $Q\leq
Z_{\frak{U}\Phi}(G)$. Thus $E$ is supersoluble by
\cite[Theorem C]{Guo181}.

\

{\bf Proof of Theorem A.}  Firstly, suppose that $\tau$ is regular. Then $X\leq
Z_{\frak{F}}(G)$ by Corollary  3.7.  Hence $E\leq Z_{\frak{F}}(G)$  by
\cite[Theorem B]{SkJGT}. Since $G/E\in \frak{F}$, we obtain $G\in \frak{F}$. Therefore, we only need to prove $G\in \frak{F}$ in the case when $\tau$ is $\Phi$-regular hereditary.

Assume
that this is false and let $(G,E)$ be a counterexample with
$|G|+|E|$ minimal. Let $F=F(E)$ and $ F^*=F^* (E)$. Let $p$ be prime
divisor of $|F^{*}|$ and   $P$ the Sylow $p$-subgroup of $F^{*}$.

(1) {\sl $X$ is  supersoluble } (This follows from Lemma 3.9).

(2) $X=F^{*}\ne E$.

Indeed, suppose that $X=E$. Then $E$ is $q$-nilpotent, where $q$ be
smallest prime divisor of $|E|$, by (1). Let $V$ be the Hall
$q'$-subgroup of $X$. If $V=1$, then $E\leq Z_{\frak{U}}(G)$ by
 Lemma 3.5, so $G\in \frak{F}$. But this contradicts the choice of $(G,
E)$. Hence $V\ne 1$. Since $V$ is characteristic in $X$, it is
normal in $G$. Moreover, the hypothesis holds for $(G/V, X/V)$ by
Lemma 2.2. Hence $G/V\in \frak{F}$  by the choice of $(G, E)$. Now we
see that the hypothesis holds also for $(G, V)$ and so $G\in
\frak{F}$ again by the choice of $(G, E)$. This contradiction shows
that we have (2).

(3) {\sl $F^{*}= F$ and  $C_{G}(F)=C_{G}(F^{*}) \leq F$}.

Since  by (1) and (2), $X=F^{*}$ is soluble,  $F^{*}=F$ by \cite[X,
 13.6]{HupBl}. We have also   $C_{G}(F)=C_{G}(F^{*}) \leq F$  by \cite[X,
 13.12]{HupBl}.

(4)  {\sl Every proper normal subgroup $W$ of $G$ with $F\leq W\leq E$ is
supersoluble}.

By   \cite[X, 13. 11]{HupBl}, $F^*(E)=F^*(F^*(E))\leq F^{*}(W)\leq F^*(E)$. It follows that
$F^{*}(W)=F^*(E)=F^*$. Thus the hypothesis is still true for $(W,
W)$ by Lemma 2.2(4). The minimal choice of $G$ implies that $W$ is
supersoluble.

(5) {\sl If $E\ne G$, then $E$ is supersoluble} (It follows directly from
(4)).

(6) {\sl If $L$ is a minimal normal subgroup of $G$ and $L\leq P$,
then $|L|>p.$}

Assume that $|L|=p$.  Let $C_0=C_E (L)$. Then the hypothesis is true
for $(G/L, C_0/L)$. Indeed, clearly, $G/C_0 =G/(E\cap C_{G}(L))$ is
supersoluble. Besides, since $L\leq Z(C_0 )$ and evidently
$F=F^{*}\leq C_{0}$  and $L\leq Z(F)$, we have
$F^{*}(C_0/L)=F^*(C_0)/L=F^{*}/L$. Hence  the hypothesis is still
true for $(G/L, C_{0}/L)$. This implies that $G/L\in \cal F$ and
thereby $G\in \cal F$ since $|L|=p$ and $\frak{U}\subseteq
\frak{F}$, a contradiction.

(7) {\sl $\Phi(G) \cap P\ne 1$ and $F^{*}(E/L)\ne F^{*}/L$ for every
minimal normal subgroup $L$ of $G$ contained in $\Phi(G) \cap P$.}

Suppose that   $\Phi(G) \cap P= 1$.  Then $P$ is the direct product
of some minimal normal subgroups of $G$ by \cite[Lemma 2.15]{Guo181}.
Hence by Lemma 3.8, $P$ has a maximal subgroup $M$ which is normal in
$G$. Now by  \cite[Lemma 2.2]{shem-skJA2009} , $G$ has a  minimal
normal subgroup $L$ with order $p$ contained in $P$, which
contradicts (6). Thus $\Phi(G) \cap P\ne 1$. Let $L\leq \Phi(G) \cap
P $ and $L$ be a minimal normal subgroup of $G$. Assume that
$F^{*}(E/L)= F^{*}/L$. Then the hypothesis is true for $G/L$ and so
$G/L\in \cal F$ by the choice of $G$. But then $G\in \cal F$ since
$L\leq \Phi(G)$. This contradiction shows that $F^{*}(E/L)\ne
F^{*}/L$.

(8)  {\sl $E$ is not soluble and $E=G$. }

Assume that $E$ is soluble.  Let $L$ be a minimal normal subgroup of
$G$ contained in $\Phi(G) \cap P$. By \cite[A, 9.3(c)]{DH},
$F/L=F(E/L)$. On the other hand, $F^{*}(E/L)=F(E/L)$ by  \cite[X,
 13.6]{HupBl}. Hence
$F^{*}(E/L)=F(E/L)= F^{*}/L$ by (3), which contradicts (7).
Therefore $E$ is not soluble and so $E=G$ by (5).

(9) {\sl   $G$ has a unique maximal normal subgroup $M$ containing
$F$, $ M$ is supersoluble and  $G/M$ is  a non-abelian simple group}
(This directly follows from (4) and (8)).

(10) {\sl  $G/F$ is  a non-abelian simple group and $G/L$ is a
quasinilpotent group if $L$ is a minimal normal subgroup of $G$
contained in $\Phi(G) \cap P$. }

Let $L$ be a minimal normal subgroup of $G$ contained in  $\Phi(G)
\cap P$.  Then by  (7), $F^{*}(E/L)\ne F^{*}/L$. Thus $F/L =F^{*}/L$
is a proper subgroup of $F^{*}(G/L)$.  By  \cite[X,
 13.6]{HupBl},
$F^{*}(G/L)=F(G/L)E(G/L)$, where $E(G/L)$ is the layer of $G/L$. By
(9), every chief series of $G$ has only one non-abelian factor. But
since $E(G/L)/Z(E(G/L))$ is a direct product of simple non-abelian
groups, we see that $F^{*}(G/L)=G/L$ is a quasinilpotent group.
Since $F(G/L)\cap E(G/L)=Z(E(G/L))$ by  \cite[X,
 (13.18)]{HupBl}, $G/F\simeq
(G/L)/(F/L)$ is a simple non-abelian group.

(11)  $F^{*}=P$.

Assume that $P\ne F$ and let $Q$ be a Sylow $q$-subgroup of $F$,
where $q\ne p$. By (10), $Q\leq Z_{\infty}(G)$. Hence by  \cite[X,
 13.6]{HupBl},
$F^{*}(G/Q)=F^{*}/Q$ and so the hypothesis is still true for $(G/Q,
G/Q)$. Hence $G/Q$ is supersoluble by the choice of $G$. It follows
that $G$ is soluble, which contradicts (8).

(12) {\sl  $p$ is the largest prime dividing $|G|$ and every Sylow
$q$-subgroup  $Q$ of $G$ where $q\ne p$ is abelian.}

Let $D=PQ$. Then $D<G$ by (8). By Lemma 3.9, $D=PQ$ is supersoluble.
Since $O_{q}(D)\leq   C_{G}(P)$,  we have $C_{G}(P)\leq P$ by (3)
and (11). Hence $O_{q}(D)=1$. Consequently, $p> q$ and $F(D)=P$.
Hence $p$ is the largest prime dividing $|G|$ and $D/P\simeq Q$ is
abelian.

{\sl Final contradiction}.

By (8) and the  Feit-Thompson theorem, 2 divides $| |G|$.  By
(12), a Sylow $2$-subgroup of $G/P$ is abelian. Hence by \cite[XI,
 13.7]{HupBl}, $G/P$ is isomorphic to one of the following: a)
$PSL(2, 2^{f})$;  b) $PSL(2, q)$, where 8 divides $q-3$ or $q-5$; c)
The Janko group $J_{1}$; d) A Ree group. It is not difficult to show
that in any case $G/P$ has a non-abelian supersoluble subgroup $V/P$
such that $p$ does not divide $|V/P|$.  Hence in view of (3) and
(11), we have $C_{V}(P)\leq P$ and so $P=F(V)$. On the other hand,
$V$ is supersoluble by Lemma 3.9. Thus $V/P$ is abelian, a
contradiction. Hence  $G\in \frak{F}$. The theorem is thus proved.

\section{Proof of Theorems B}

{\bf  Lemma 4.1.} {\sl  Let $\frak{F}$ be a saturated formation, $P$
be a normal $p$-subgroup of a group $G$, where $p$ is a prime. Let
$D$ be a characteristic subgroup of $P$ such that every non-trivial
$p'$-automorphism of $P$ induces a non-trivial automorphism of $D$.
Suppose that $D \leq Z_{{\cal F}}(G)$. Then $P \leq
Z_{\frak{F}}(G)$.  }

{\bf Proof.} Let $F$ is the canonical local satellite of $\frak{F}$.
Let $C=C_{G}(P)$.  Since $D \leq Z_{\frak{F}}(G)$, then
$G/C_{G}(D)\in F(p)$ by Lemma 3.2. On the other hand, since every
non-trivial $p'$-automorphism of $P$ induces a non-trivial
automorphism of $D$,  $C_{G}(D)/C_G(P)$ is a $p$-group. Hence from
the definition of $F$ we have  $G/C_{G}(P)\in F(p)$. It follows that $P \leq
Z_{\frak{F}}(G)$.

Let $P$ be a non-identity  $p$-group. If $P$ is not a non-abelian
2-group, we use $\Omega (P)$ to denote the subgroup $\Omega_{1} (P)$.
Otherwise, $\Omega (P)=\Omega_{2} (P)$.

{\bf  Lemma 4.2} (see \cite[Theorem 2.4]{Gag}). {\sl Let $P$ be a
group, $\alpha $ a $p'$-automorphism of $P$. If $[\alpha,
\Omega(P)]=1$, then $\alpha =1$. }

{\bf  Lemma 4.3.} {\sl  Let $C$ be a Thompson  critical subgroup of
a $p$-group $P$} (see \cite[p. 185]{Gor}). {\sl Then the group
$D:={\Omega}(C)$ is of exponent $p$ if $p$ is an odd prime, or
exponent 4 if $P$ is non-abelian 2-group. Moreover, every
non-trivial $p'$-automorphism of $P$ induces a non-trivial
automorphism  of $D$.}

{\bf Proof.}  See the proof of \cite[Chapter 5, 3.13]{Gor} and use the
fact that if $C$ is a non-abelian 2-group, then ${\Omega}(C)$ is of
exponent 4 (see \cite[p.~3]{rerk-kas}).

{\bf  Lemma 4.4.} {\sl Let $\frak{F}$ be a saturated formation, Let
$P$ be a normal $p$-subgroup of a group $G$ and $D={\Omega}(C)$,
where $C$ is  a Thompson critical subgroup of $P$. If  $C \leq
Z_{{\cal F}}(G)$  or $D \leq Z_{\frak{F}}(G)$, then $P \leq
Z_{\frak{F}}(G)$. }

{\bf Proof.} Let $Z=Z_{\frak{F}}(G)$. Suppose  that  $C \leq Z$.
Then $G/C_{G}(C)\in F(p)$, where  $F$ is the canonical local
satellite of $\frak{F}$ by  Lemma 3.2. On the other hand,
$C_{G}(C)/C_{G}(P)$ is a $p$-group by \cite[Chapter 5, 3.11]{Gor}. Hence
$G/C_{G}(P)\in F(p)$. Consequently, $P \leq Z$. On the other hand,
by Lemmas 4.2 and 4.3, $C_{G}(C)/C_{G}(D)$ is also a $p$-group, and so in the case
when $D\leq Z$ we similarly get that $C\leq Z$.

{\bf  Lemma 4.5.} {\sl  Let $P/R$ be a chief factor of a group $G$
with $|P/R|= p^{n}$, where  $p$ is a prime and $n > 1$. Suppose that
for every normal subgroup $V$ of $G$ with $V < P$ we have $V\leq R$.
Let $H$ be a subgroup of $P$ such that $R < RH < P$. If  $H$ is a
cyclic group of prime order or order 4, and $T$ is a supplement of
$H$ in $G$, then $T=G$.}

 {\bf Proof.}   Assume that $T\ne G$. Then
$P=H(P\cap T)$, where $P\cap T\ne P$ and $|P:P\cap T|=|G:T|$. Let
$N=N_{G}(P\cap T)$. Since $T\leq N$ and $P\cap T <  N_{P}(P\cap T)$,
we have that either $P\cap T$ is a normal subgroup  of $G$  or
$|G:N|=2$. The first case implies that $P\cap T\leq R$ and hence
$P=RH$, a contradiction. In the second case, $N$ is normal in $G$
and  so $N\cap P$ is a normal subgroup of $G$ with $|P:P\cap N|=2.$
Therefore $P\cap N\leq R$ and thereby $|P/R|=2$, a contradiction.
Hence $T=G$.

{\bf  Proposition 4.6.}  {\sl Let $\frak{F}$ be a  saturated
formation containing all supersoluble groups and  $P$  a
non-identity normal $p$-subgroup of   $G$. Suppose that
every cyclic subgroup of $P$ of prime order or order 4 (if $P$ is a non-abelian group)
is  ${\frak{F}}_{\tau}$-supplemented in $G$.}

(a) {\sl  If $\tau$ is $\Phi$-quasiregular and  $P$ is of exponent $p$ or  exponent $4$, then $P\leq
Z_{\frak{F}\Phi}(G)$.}

(b) {\sl If $\tau$ is quasiregular, then $P\leq Z_{\frak{F}}(G)$.}

{\bf Proof.}  Suppose that in this  theorem is false and let $(G,
P)$ be a counterexample with $|G|+|P|$ minimal. We write
$Z=Z_{\frak{F}\Phi}(G)$ if $\tau$ is $\Phi$-quasiregular and $P$ is of exponent $p$ or  exponent $4$,
and  $Z=Z_{\frak{F}}(G)$ if $\tau$ is quasiregular. Let   $G_{p}$ a
Sylow $p$-subgroup of $G$.

(1) {\sl $G$  has a  normal subgroup $R\leq P$ such that $P/R$ is an
$\frak{F}$-eccentric chief factor of $G$, $R\leq Z$  and $V\leq R$
for any normal subgroup $V\ne P$ of  $G$  contained in $P$. In
particular, $|P/R| > p$.  }

Let $P/R$ be a chief factor  of $G$. Then  the hypothesis holds for
$(G, R)$. Therefore $R\leq Z$  and so $P/R$ is
$\frak{F}$-eccentric in $G$ by the choice of $(G, P)$ and \cite[Lemma
2.2]{shem-skJA2009}. It follows that $|P/R|>p$. Now let $V\ne P$
be any normal subgroup of $G$ contained in $P$. Then $V\leq Z$. If
$V\nleq R$, then $VR=P\leq Z$ by \cite[Lemma
2.2]{shem-skJA2009}.   This
contradiction shows that $V\leq R$.

(2)  {\sl  $P$ is of exponent $p$ or  exponent $4$.}

Assume that this is false and let $\tau$ be quasiregular. Let $L$ be a Thompson critical subgroup of $P$ and
$\Omega =\Omega (L)$. Then $\Omega < P$, and so $ \Omega
\leq Z$ by (1), where $Z=Z_{\frak{F}}(G)$. Hence $P\leq Z$ by Lemma 4.4,
which contradicts the choice of $(G, P)$. Hence $\Omega  = P$, so we
have (2) by Lemma 4.3.

Now let $L/R$ be any minimal subgroup of $(P/R)\cap Z(G_{p}/R)$. Let
$x\in L\setminus R$ and $H=\langle x \rangle$. Then $L/R=HR/R$, so
$H$ is not normal in $G$ and $H_{G}\leq R$ by (1). Moreover, $|H|$
is ether prime or 4 by (2). Hence  $H$ is
${\frak{F}}_{\tau}$-supplemented in $G$, so there is a subgroup
 $T/H_{G}$ of $G/H_{G}$ such that
$(T/H_{G})(H/H_{G})=G/H_{G}$ and $(T/H_{G})\cap (H/H_{G}) \subseteq
(S/H_{G})Z_{\cal F}(T/H_{G})$  for some $\tau$-subgroup $S/H_{G}$ of
$G/H_{G}$   contained in $H/H_{G}$.

(3) {\sl $H/H_{G}=S/H_{G}$  and $L/R\in \tau (G/R)$ }.

By Lemma 4.5, we have $T=G$. Hence
$H/H_{G}=(S/H_{G})(H/H_{G}\cap Z_{\frak{F}}(G/H_{G}))$. On the
other hand, since $H$ is cyclic, we have either $H/H_{G}\leq
Z_{\frak{F}}(G/H_{G})$ or $H/H_{G}=S/H_{G}$. Note that
$(R/H_{G})(H/H_{G})/(R/H_{G})=(RH/H_{G})/(R/H_{G})=(L/H_{G})/(R/H_{G})$,
so if we have the former case, then $(L/H_{G})/(R/H_{G})\leq
(P/H_{G})/(R/H_{G})\cap Z_{\frak{F}}((G/H_{G})/(R/H_{G}))$ by
\cite[Lemma 2.2(4)]{GuoS-JA}. But then  $L/R\leq (P/R)\cap
Z_{\frak{F}}(G/R)$. This implies that $P/R$ is ${\frak{F}}$-central in $G$,
contrary to (1). Therefore we have $H/H_{G}=S/H_{G}$, so
$(L/H_{G})/(R/H_{G})=(H/H_{G})(R/H_{G})/(R/H_{G})$ is a
$\tau$-subgroup of $(G/H_{G})/(R/H_{G})$ and hence $L/R$ is a
$\tau$-subgroup of $G/R$.
 Hence we have (3).

(4) {\sl $\tau$ is not quasiregular.}

Assume that this is false. In view of (3), $L/R$ is a
$\tau$-subgroup of $G/R$. But since $L/R$ is normal in $G_{p}/R$,  $L/R$ is normal in $G/R$ by Lemma 3.4. Hence $L/R=P/R$,
which contradicts (1). Hence we have (4).

Now, in view of (4), we have only to prove that $P\leq
Z_{\frak{F}\Phi}(G)$.

(5) {\sl There is a maximal subgroup $M$ of $G$ such that $R=P\cap
M$ and $MP=G$}.

Indeed, if $P/R\leq \Phi (G/R)$, then in view of  (1) and
\cite[Lemma 2.2]{GuoS-JA}, $P\leq Z_{\frak{F}\Phi}(G)$, which
contradicts the choice of $(G, P)$. Therefore for    some maximal
subgroup $M/R$ of $G/R$ we have $(M/R)(P/R)=G/R$. Then $PM=G$ and
$M\cap P=R$ since $P/R$ abelian.

{\sl Final contradiction.}

Let $D=M_{G}$. Clearly, $R\leq D$. By (3), $LD/D$ is a $\tau$-subgroup  of $G/D$. On
the other hand, since $L\nleq R$,  $L\nleq D$ and so $|LD/D|=|L/L\cap
D|=|L/R|=p$. Since $\tau$ is $\Phi$-qusiregular,  $|G/D: N_{G/D}( LD/D)|$
is a power of $p$. On the other hand, since $L/R$ is normal in
$G_{p}/R$, $LD/D$ is normal in $G_{p}D/D$. Hence in view of Lemma 3.4,
$L^{G}D=LD$. This induces that $PD=LD$ and $P=L(P\cap D)=LR=L $, which contradicts (1).
The proposition is thus proved.

{\bf  Lemma 4.7}   {\sl If  $G=NT$, where $T$ is a proper
subgroup of $G$, $N\leq Z_{\cal \infty}(G)$ and $N$ is a $p$-group,
then $O^{p}(G)\ne G$.}

{\bf Proof.} Let $T\leq M$ where $M$ is a maximal subgroup of $G$.
Then $G/M_{G}=(NM_{G}/M_{G})(M/M_{G})$ is a primitive group and $NM_{G}/M_{G}\leq Z_{\cal
\infty}(G/M_{G})$. Hence $NM_{G}/M_{G}\leq Z(G/M_{G})$
and so $M/M_{G}$ is normal in $G/M_{G}$. Therefore $O^{p}(G)\leq M$.

{\bf Lemma 4.8.}   {\sl Let $G=NT$, where $N$ is a minimal normal
subgroup of $G$ and  $T$ is a maximal subgroup of $G$.}

(1) {\sl If $|G:T|$ divides 4, then $N$ is abelian}.

(2) {\sl If $N\leq Z_{\cal U}(G)$, then $|G:T|$  is a prime.}

{\bf Proof.} (1) If $|G:T|=4$, then $G/T_{G}$ is isomorphic with a
subgroup of the symmetric group $S_{4}$ of degree 4, which implies
that $N\simeq NT_{G}/T_{G}$ is abelian. If $|G:T|=2$, then clearly
$N\leq Z(G)$.

(2) Assume that $N\leq Z_{\cal U}(G)$. Since $N$ is a minimal normal
of $G$, the order of $N$ is a prime. Thus $|G:T|$  is a prime.

{\bf  Theorem 4.9.}  {\sl Let $\frak{F}$ be the class of all
$p$-supersoluble groups. Let     $E$ be  a normal subgroup
$G$ and $P$ a Sylow  $p$-subgroup of $E$ of order $|P|=p^{n}$, where $n > 1$
and $(|E|, p-1)=1$. Suppose that every cyclic subgroup of $P$ of
prime  order or 4 (if $P$ is a non-abelian 2-group) is
${\frak{F}}_{\tau}$-supplemented in $G$. If $\tau$ is either
hereditary or regular, then $E$ is $p$-nilpotent.}

{\bf Proof.} Suppose that this theorem is false and  let $(G, E)$
be a counterexample with $|G|+|E|$  minimal. Let $Z=Z_{\cal F}(G)$.

(1) $O_{p'}(E)=1$ (see (1) in the proof of Theorem 3.6).

(2) {\sl $\tau$ is regular.}

Assume that this is false. Then, by hypothesis, $\tau$ is
hereditary. Therefore the hypothesis holds for every subgroup $B$ of
$G$ containing $P$. Hence $E=G$  and every maximal
subgroup of $G$ is $p$-nilpotent by the choice of $G$. Hence by
\cite[IV, 5.4]{hupp},  $G=P\rtimes Q$ is a $p$-closed Schmidt
group, where $Q$ is a Sylow $q$-subgroup of $G$ ($q\ne p$) and   $P$
is of exponent $p$ or  exponent 4 (if $P$ is a non-abelian 2-group).
But then, by Proposition 4.6, $P\leq Z_{\frak{U}\Phi}(G)=Z_{\infty}(G)$
since $(|G|, p-1)=1$. Hence  $G$ is nilpotent. This contradiction
shows that we have (2).

(3) $O_{p}(E)\leq Z_{\infty}(E)$.

Assume  $O_{p}(E)\ne 1$. Let $H$ be a subgroup of $O_{p}(E)$ which
is ${\frak{F}}_{\tau}$-supplemented in $G$. Then there is a subgroup
$T/H_{G}$ of $G/H_{G}$ such that $(H/H_{G})(T/H_{G})=G/H_{G}$ and
$(H/H_{G})\cap (T/H_{G}))\subseteq (S/H_{G})Z_{\cal F}(T/H_{G})$ for some
$\tau$-subgroup $S/H_{G}$ of $G/H_{G}$ contained in $H/H_{G}$. Then
$(H/H_{G})\cap (T/H_{G}))\subseteq (S/H_{G})(Z_{\cal F}(T/H_{G})\cap
(O_{p}(E)/H_{G}))\leq (S/H_{G})Z_{\cal U}(T/H_{G})$. Hence $H$ is
${\frak{U}}_{\tau}$-supplemented in $G$. Now since $\tau$ is
regular, $O_{p}(E)\leq Z_{\cal U}(G)$ by Lemma 4.6. But since $(|E|, p-1)=1$, $Z_{\frak{U}}(G)=Z_{\infty}(G)$. Thus $O_p(E)\leq Z_{\infty}(G)\cap E\leq Z_{\infty}(E)$.

(4) {\sl If $V/D$ is a chief factor of $G$ where $D\leq O_{p}(E)$,
 then $V\leq O_{p}(E)$.}

Assume that this is false and let $(V, D)$ be the pair with
$|V|+|D|$ minimal such that $V/D$ is a chief factor of $G$, $D\leq
O_{p}(E)$ and $V\nleq O_{p}(E)$. Then:

(a) {\sl $p$ divides $|V/D|$.}

Assume that $V/D$ is a $p'$-group. Then  $V$ is $p$-soluble. By (3),
$D\leq Z_{\infty}(E)\cap V\leq Z_{\infty}(V)$. Hence $V/C_{V}(D)$ is
a $p$-group by Lemma 3.2. Hence for a
Hall $p'$-subgroup $W$ of $V$ we have $V=D\times W$. Then  $W$ is
characteristic in $V$. Hence $W\leq O_{p'}(E)$, contrary to (1).
Therefore we have (a).

(b) {\sl $p=2$,  $V/D$ is  non-abelian and $O_{2}(E)\leq
Z_{\infty}(G)$.}

Since $V\nleq O_{p}(E)$, the choice of $(V, D)$ and Claim (a) imply
that $V/D$ is  non-abelian. But since $(|E|, p-1)=1=(|V|, p-1)$,  by the
Feit-Thompson's theorem we have $p=2$. By Lemma 4.6, $O_{2}(E)\leq
Z_{\infty}(G)$.

(c) {\sl $V$ has a 2-closed Schmidt subgroup $A$ such that for some
cyclic  subgroup $H$ of $A$ of order 2 or order 4 we have $H\nleq
D$. }

In view of (b) and \cite[Theorem 3.4.2]{Guo}, $V$ has a 2-closed Schmidt
subgroup $A$ and for the  Sylow 2-subgroup $A_{2}$ of $A$ the following
hold: (i) $A_{2}=A^{\frak{N}}$; (ii) $A_{2}$ is of exponent 2 or
exponent 4 (if $A_{2}$  is non-abelian); (iii) $\Phi
(A)=Z_{\infty}(A)$; $A_{2}/\Phi (A_{2})$ is a non-cyclic chief
factor of $A$. Therefore for some cyclic  subgroup $H$ of $A$ of
order 2 or order 4 we have $H\nleq \Phi (A_{2})$. But $\Phi
(A_{2})=Z_{\infty}(A)\cap A_{2}$, so  $H\nleq Z_{\infty}(V)\cap
A\leq Z_{\infty}(A)$, which implies that $H\nleq D$ by (3). It is
also clear that $H_{G}\leq D$.

Without loss of generality, we may assume that $H\leq P$, so $H$ is
${\frak{F}}_{\tau}$-supplemented in $G $.

For any subgroup $L$   of $G$, we  put $\widehat{L}=LD/D$.

(d) {\sl For any non-identity subgroup $S/H_{G}$  of $G/H_{G}$ not
contained in $D/H_{G}$ we have $S/H_{G}\not \in \tau(G/H_{G})$.}

Suppose that this is false. Sincee $\tau $ is regular,
$|\widehat{G}:N_{\widehat{G}}(\widehat{H})|$ is a power of 2, and so
$\widehat{H}\leq O_{2}(\widehat{G})$ by Lemma 3.4. Hence $\widehat{V}\leq
O_{2}(\widehat{G})$, a contradiction. Thus (d) holds.

(e) {\sl For any $\tau$-subgroup $S/H_{G}$ of $G/H_{G}$ contained in
$H/H_{G}$ we have
 $H/H_{G}\nleq    (S/H_{G})Z_{\frak{F}}(G/H_{G}).$ }

Suppose that $H/H_{G}\leq    (S/H_{G})Z_{\frak{F}}(G/H_{G})$. Then
$H/H_{G}=(S/H_{G})((H/H_{G}) \cap Z_{\frak{F}}(G/H_{G}))$. Since
$H/H_{G}$ is cyclic, either $H/H_{G}=S/H_{G}$ or $H/H_{G}\leq
Z_{\frak{F}}(G/H_{G})$. But the former case  is impossible by
(d). Hence  $H/H_{G}\leq Z_{\frak{F}}(G/H_{G})$. But then
$(H/H_{G})(D/H_{G})(D/H_{G})\leq Z_{\frak{F}}((G/H_{G})/(D/H_{G}))$
by \cite[Lemma 2.2(4)]{GuoS-JA}. Hence $V/D\simeq (V/H_{G})/(D/H_{G})$
is a 2-group, a contradiction. Therefore we have (e).

(f) {\sl $O^{2}(V)= V$. }

Assume that  $O^{2}(V)\ne V$. Since $O^{2}(V)$ is characteristic in
$V$, it is normal in $G$. Hence $DO^{2}(V)=V$. Then in view of the
$G$-isomorphism $V/D\simeq O^{2}(V)/D\cap O^{2}(V)$, $O^{2}(V)/D\cap
O^{2}(V)$ is a non-abelian chief factor of $G$ such that
$O^{2}(V)\nleq O_{2}(V)$ and  $D\cap O^{2}(V) \leq O_{2}(E)$, which
contradicts the choice of $(V, D)$. Hence we have (f).

{\sl final contradiction.}

Since  $H$ is ${\frak{F}}_{\tau}$-supplemented  in $G $, there is a
subgroup  $T/H_{G}$  of $G/H_{G}$ such that
$(H/H_{G})(T/H_{G})=G/H_{G}$ and $(H/H_{G})\cap (T/H_{G})\subseteq
(S/H_{G})Z_{\frak{F}}(T/H_{G})$ for some $\tau$-subgroup $S/H_{G}$
of $G/H_{G}$ contained in $H/H_{G}$.

Suppose that $T=G$. Then $(H/H_{G})\leq
(S/H_{G})Z_{\frak{F}}(G/H_{G})$,  which contradicts (e).
Therefore $T\ne G$.

Since $HT=G$,  $|G:T|$ divides 4. Let $T\leq M$, where $M$ is a
maximal subgroup of $G$. Then $V\nleq M$. Suppose that $|G:M|=2$.
Then from the isomorphism $G/M\simeq V/V\cap M$ we get $O^{2}(V)\ne
V$, which contradicts (f). Hence $M=T$ is a maximal subgroup of $G$ and
$|G:T|=4$.

Note that if $TD=G$, then $V=D(T\cap V)$. Hence $O^{2}(V)\ne V$ by
Lemma 4.7, contrary to (f). Therefore $TD\ne G$, which in view of
maximality of $T$ implies that $D\leq T$. Therefore
$\widehat{G}=\widehat{H}\widehat{T}=\widehat{V}\widehat{T}$. Then by Lemma 4.8,
$\widehat{V}$ is abelian, which contradicts (b). This final contradiction completes the proof.

{\bf  Corollary 4.10.}  {\sl  Let $E$ be a non-identity normal
subgroup $G$. Suppose that for every non-cyclic Sylow subgroup $P$ of $E$, every cyclic subgroup of $P$ of
prime  order or 4 (if $P$ is a non-abelian 2-group) is
${\frak{U}}_{\tau}$-supplemented in $G$ for sone hereditary
quasiregular subgroup functor  $\tau $. Then $E\leq
Z_{\frak{F}}(G)$.}

{\bf Proof.}  See the proof of Corollary 3.7 and use  Proposition 4.6 and Theorem 4.9
instead of  Proposition 3.5 and Theorem 3.6,  respectively.

{\bf Proof of Theorem B.} Suppose that this theorem is false and let
$(G, E)$ be a counterexample with $|G|+|E|$ minimal. Then
$E=G^{\frak{F}}$. Let $p$ be the smallest
prime dividing $|X|$. Then $X$ is $p$-nilpotent by Theorem 4.9 and
\cite[Chapter 7, 6.1]{Gor}. Hence $X$ is soluble by the Feit-Thompson's
theorem. Let $V$ be the Hall $p'$-subgroup of $X$. Then $V$ is
characteristic in $X$ and so it is normal in $G$.

Suppose that Assertion (i) is true. Then the choice of $G$ implies
that $E\nleq Z_\frak{F}(G)$ and so $V\ne
1$ by Proposition 4.6. It is clear that the hypothesis holds for $(G/V, E/V)$ by Lemma 2.2. Hence $G/V\in \frak{F}$ by the choice of $(G, E)$. Then the hypothesis holds
for $(G, V)$. Hence the choice of $(G, E)$ implies that $G\in
\frak{F}$, a contradiction. Hence at least one of the Assertions (ii)
or (iii) is true. In this case $\tau$ is either hereditary
quasiregular or regular. By Corollary 4.10, $X\leq Z_{\frak{U}}(G)$.
It follows from \cite[Proposition C]{SkJGT} that $E\leq Z_{\frak{U}}(G)$. Then
$G\in \frak{F}$. This contradiction completes the proof.

\section{ Final remarks}

I. A large number of known results are corollaries of Theorems A and B. In particular, in view of Example 1.4, Theorem A covers Theorem D in \cite{Ezq1} and Theorem 4.7 in \cite{Weihua1};
 in view of Example 1.5, Theorems B covers Theorem 1.3
in \cite{Guo8};
in view of the remarks after Definition 1.2, Theorems A and B
cover Theorems 1.2 and 1.4 in \cite{ahmad},  Theorem 3.1 and 3.2 in
\cite{GuoManu1} and Theorem 3.3 in \cite{Acta};
in view  of Example 1.6 and the  remarks after Definition 1.2,
Theorems A and B cover Theorem 4.1 in \cite{Ballest}, Theorems 1.1 and 1.2 in \cite{Weihua}, Theorems 3.1,
3.4 and 3.6 in \cite{long}, Theorem A in \cite{SkJGT} and Theorems A
and B in \cite{China};
in view of Example 1.7, Theorems A and B cover Theorems 5.1 and 5.2 in
\cite{GuoS176};
in view of Example 1.8, Theorem B covers Theorems 3.1 and 3.2 in \cite{Lixian};
in view of Example 1.9,
 Theorem A and B cover   Theorem 1.3   in
\cite{Ukr1} and Theorem 1.2     in
\cite{Ukr2};
in view of Example 1.10,
 Theorems A and B cover   Theorems 3.1, 3.4-3.7 in
\cite{LiS} and so on.

II. We do not know now the answer to the following

{\bf  Question B.} {\sl Is  the subgroup functor in (1.5)
regular?}  (see Question 17.112 in  \cite{Kour}).

\end{document}